\def\a{\alpha}
\def\b{\beta}

\def\g{\gamma}
\def\eps{\varepsilon}

\def\O{\Omega}

\def \dd#1{{\bf#1}}

\def\cl#1{{\cal#1}}



\def\ouv#1{\smash{\mathop{#1}\limits^{\lower 1pt\hbox
{$\scriptscriptstyle\circ$}}}}

\def\apl#1#2#3#4#5{{#1}:\matrix{
		\hfill#2 & \longrightarrow & #3\hfill \cr
		\hfill#4 & \longmapsto     & #5\hfill \cr}}

\def\hfl#1#2{\smash{\mathop{\hbox to 12mm{\rightarrowfill}}
\limits^{\scriptstyle#1}_{\scriptstyle#2}}}


\long\def\eno#1#2{\par\smallskip{\bf{#1}}{\it\ {#2}}\par\medskip}

\def\tit#1{\vskip 5mm plus 1mm minus 2mm {\tir #1}
		\vskip 3mm plus 1mm minus 2mm}

\def\stit#1{\vskip 3mm plus 1mm minus 2mm {\bf{#1}}
		\smallskip}

\long\def\dem#1#2{{\bf {#1}}{\ {#2}$\diamondsuit$}\medskip}

\font\tir=cmbx10 at 12pt

\def\ref#1#2#3#4{{\bf #1}{\ #2}{\it ,\ #3}{,\ #4}\medskip}


\def \picture #1 by #2 (#3){\midinsert \centerline 
{\vbox to #2{\hrule width #1 heigth 0pt 
depth 0pt \null \vfill \special {picture #3}}}\endinsert }

\def\scaledpicture #1 by #2 (#3 scaled #4) {{
\dimen0 =#1 \dimen1 =$2
\divide \dimen0 by 1000 \multiply \dimen0 by #4
\divide \dimen1 by 1000 \multiply \dimen1 by #4
\picture \dimen0 by \dimen1 (#3 scaled $4)}}

\def\figure #1 #2 #3 {\midinsert \vglue 3mm 
{\vbox to #3 {\hrule width 6cm height 0cm depth 0cm \vfill
{\special {picture #1 scaled #2}}}}\vglue 2mm \endinsert}

\magnification=1200


\overfullrule=0pt

{\centerline {\tir {Total convergence or general divergence}}}

{\centerline {\tir {in Small Divisors}}}

\bigskip
\bigskip

{\centerline {{\bf R. P\'erez-Marco}\footnote {*} 
{  UCLA, 
Dept. of Mathematics, 
405, Hilgard Ave., Los Angeles, CA-90095-1555, 
USA, e-mail: ricardo@math.ucla.edu; CNRS UMR 8628, 
Universit\'e Paris-Sud, Math\'ematiques, 91405-Orsay, France.}}}

\bigskip


{\bf Abstract.} {\it We 
study generic holomorphic 
families of dynamical systems presenting  
problems of small divisors with
fixed arithmetic. We prove that we have convergence for all 
parameter values or divergence everywhere except for 
an exceptional set of zero $\Gamma$-capacity. 
We illustrate this general principle in  
different problems of small divisors. As an application
we obtain
new richer families of non-linearizable examples in the 
Siegel problem when Bruno condition is violated, generalizing 
previous results of Yoccoz and the author.
}

\bigskip


\medskip


\bigskip
\bigskip

\tit {1) Introduction.}

In this article we study generic (polynomial) 
holomorphic families of dynamical systems 
presenting problems of 
small divisors with
fixed arithmetic. The principle our
theorems illustrate is the following : 

\medskip

{\bf We have total convergence for all 
parameter values or general divergence except maybe for a 
very small exceptional set of parameter values.}

\medskip

The germinal idea can be traced back to Y. Ilyashenko 
where in [Il] he studies
divergence in problems of small divisors from  
divergence of the homological (or linearized) equation.  
Ilyashenko's paper contains a remarkable idea.
We find there, for the first
time in Small Divisors, 
the study of linear deformation of the system and the use
of the polynomial dependence of the new formal linearizations.
A similar idea, but not quite in the same problem,
was used by H. Poincar\'e to show that linear deformations
of completely integrable hamiltonians are not generally
completely integrable with analytic first integrals depending
analytically on the parameter ([Poi] volume I chapter V). 
It is worth noting that 
this is the key preliminary step in his difficult proof of the 
non existence of non trivial local analytic first integrals in the three 
body problem.

Such a linear deformation 
has been fruitfully used by J.-C. Yoccoz. He proves that  
in the Siegel problem
the quadratic polynomial is the worst linearizable holomorphic
germ ([Yo] p. 58).
The only ingredient in this proof that is not 
in Ilyashenko's one is the classical Douady-Hubbard
straightening theorem for polynomial-like mappings.
Yoccoz simplifies Ilyashenko's argument 
replacing Nadirashvili's lemma
by the maximum principle. He loses in that way the strength
of the original approach, in particular the potential theoretic
aspects. Non-linear polynomial perturbations 
were used by the author in [PM1] to 
generalize Yoccoz's result to higher degree polynomials.

\medskip

In this article we clarify and strength the role played by 
potential theory in  parameter space. 
A key point is the observation 
that Nadirashvili's lemma  
can be improved by using Bernstein's lemma in 
approximation theory. In that way we do make precise the
thinest notion for the exceptional set. In parameter space
the exceptional set has  $\Gamma$-capacity $0$
(for the definition of $\Gamma$-capacity see [Ro] section 2.2).
The intersection of the exceptional 
set with any complex line is full (the whole line) or a polar set
(any set of $\Gamma$-capacity $0$ has this property, 
see [Ro] Lemma 2.2.8 p.92).

\medskip

As far as the author knows, the first person who studied 
small divisors problems using ingredients from potential
theory in parameter space is M. Herman (see [He1] and [He2]).

\medskip

The techniques in this article are applicable to virtually 
any holomorphic problem in small divisors where the dependence 
on parameters of the coefficients of the divergent
series are polynomial (as we will see this happens in most of 
the problems). We 
have selected a few illustrative ones guided 
mainly by our personal taste.
We only consider here polynomial families. The same proof 
can be extended easily for more general holomorphic families 
(see remark 4).

It is surprising   
that the idea and the results of this article 
have been overlooked so far.

\stit {Linearization.}

\eno {Theorem 1.}{ Let $n,m\geq 1$ and $d\geq 0$. For a multi-index
$i=(i_1, \ldots
i_m) \in \dd N^m$ with $0\leq i_1+\ldots i_m \leq d$, let 
$f_i$ be a germ of holomorphic map 
$$
f_i : (\dd C^n , 0) \to (\dd C^n , 0)
$$
with valuation larger or equal to $2$ (i.e. $f_i (z)=\cl O (z^2)$).

For $t=(t_1, \ldots t_m) \in \dd C^m$
we consider the holomorphic family of germs of holomorphic 
diffeomorphisms, $z\in \dd C^n$,
$$
f_t(z)= Az +\sum_{i=(i_1, \ldots i_m)\atop 0\leq i_1+\ldots +i_m \leq d} 
t^i f_i (z)
$$
where $A\in GL_n(\dd C )$ is a fixed linear map, $A=D_0 f$, 
with non-resonant eigenvalues.

Then all maps $f_t$, $t\in \dd C^m$ are formally 
linearizable, i.e. there exists a unique formal map $h_t$ with $h_t(0)=0$
and $D_0h_t=I$ such that the formal equation 
$$
h_t\circ f_t  =A\circ h_t 
$$
is satisfied.

We have the following dichotomy:

\medskip

1) The holomorphic family $(f_t)_{t\in \dd C^m}$ 
is holomorphically linearizable,
that is for all $t\in \dd C^m$, $h_t$ defines a germ of holomorphic 
diffeomorphism. Moreover, the radius $R(t)$ of convergence of the 
linearization $h_t$ is bounded from below on compact sets 
and, more precisely,  for some $C_0 >0$,
$$
R(t) \geq {C_0 \over 1+||t||} \ .
$$

\medskip

2) Except for an exceptional set  
$E\subset \dd C^m$ of $\Gamma$-capacity $0$ of values of $t$, 
$f_t$ not holomophically linearizable.

}

{\bf Remarks:}

\medskip

1) We remind that the eigenvalues 
$(\lambda_1 , \lambda_2, \ldots , \lambda_n)$
of $A$, counted with multiplicity, are non-resonant if 
$$
\lambda_i-\lambda_1^{i_1}\ldots \lambda_n^{i_n} \not=0
$$
for all $(i_1, \ldots , i_n)\in \dd N^n$ with $i_1+\ldots +i_n \geq 2$.
We give later a theorem for holomorphic germs with resonant linear parts.

\medskip

2) The linear part $A$ is in the Poincar\'e domain if 
$$
\min (\max_i |\lambda_i| ,\max_i |\lambda_i^{-1}| ) <1 \ .
$$
In that case it is well known that we are always in case (1).
Otherwise the linear part of $A$ belongs to the Siegel domain.

\medskip

3) In general the exceptional set 
$E\subset \dd C^m$ is not empty. For
example if $f_0=0$ and $0\in W$, then $0\in E$ when we are 
in the second case. 

\medskip

4) With the same type of proof, one can prove the result for holomorphic
families of the form
$$
f_t(z)= Az +\sum_{i=(i_1, \ldots i_m)\atop 0\leq i_1+\ldots +i_m } 
t^i f_i (z)
$$
where the holomorphic germs $(f_i(z))$ have valuations such that 
${\hbox {\rm val}}f_i \geq \eps_0 |i|$ for some $\eps_0 >0$.

\medskip

Some illustrative corollaries follow now. 
For $n=1$ and the special case of entire 
functions we have directly from theorem 1:

\eno {Corollary 1.}{Let $(f_t)_{t\in \dd C^m}$ be a finite
dimensional holomorphic family of entire functions as above with
$$
f'_t(0)=\lambda 
$$
where $\lambda=e^{2\pi i \a }$ with $\a \in \dd R-\dd Q$.

Then the family is linearizable or, except for an exceptional polar 
set $E\in \dd C$ of values of $t$, all $f_t$ are non-linearizable.}

Assuming that the family contains a non-linearizable structurally 
stable polynomial (for example a quadratic polynomial) we can break 
the dichotomy. This just follows from the observation that in a 
neighborhood of 
this polynomial all elements of the family are quasi-conformally 
conjugated (by Douady-Hubbard straighthening theorem), thus 
they are linearizable or not simultaneously. 

\eno {Corollary 2.}{Let Let $(f_t)_{t\in \dd C^m}$ be a finite
dimensional holomorphic family of entire functions as above with
$$
f'_t(0)=\lambda 
$$
where $\lambda=e^{2\pi i \a }$ with $\a \in \dd R-\dd Q$.
We assume that for some value $t_0$ $f_{t_0}$ is a structurally 
stable polynomial in the space of polynomials with fixed point 
at $0$ and multiplier $\lambda$.

Then if $\a$ is not a Bruno number almost all entire functions 
$f_t$, except maybe for an exceptional polar set $E\subset \dd C$ of 
values of $t$,  
are not linearizable.}

When $\a \in \dd R-\dd Q$ is not a Bruno number, no examples 
were known of non-linearizable entire 
functions not quasi-conformally conjugated to polynomials in 
a neighborhood of $0$. This was due to the shortcomings of 
Yoccoz maximum principle approach [Yo].

A particular case of this corollary is  the theorem
proved in [PM1] about polynomial germs. The author showed,
generalizing Yoccoz result for the quadratic polynomial, 
that if $\a$ is not a Bruno number,
in the space
$$
\cl P_{\lambda, d}=\{ P(z)=\lambda z+a_2 z^2 +\ldots +a_d z^d ; 
(a_2, \ldots , a_d ) \in \dd C^{d-1} \} 
$$
the polynomials that are not structurally stable (this is an
open dense set whose complement has $\Gamma$-capacity $0$) 
are not linearizable. 

It is worth mentioning that the question to know if the 
exceptional set $E_{\lambda ,d}$ is empty for the polynomial family
$\cl P_{\lambda , d}$ when $\a \in \dd R -\dd Q$ is not a 
Bruno number,
is still open, even for the cubic family:
$$
P_b (z)=\lambda z +bz^2 +z^3 \ 
$$
Contrary to unanimous belief, the author will not be surprised
that $E_{\lambda , d}$ is not empty for appropriate values of 
$\lambda$ and $d$. For Liouville numbers $\a$ with extremely 
good rational approximations, by an argument
of Cremer (see [Cr] and [PM1]), $E_{\lambda , d}$ is known to 
be empty.

To illustrate the strength of the precedent theorem, we present
the following variations.

\eno {Corollary 3.}{Let $\a \in \dd R-\dd Q$ be not Bruno.

1) Let $f(z)=e^{2\pi i \a} z +\cl O (z^2 )$ be non-linearizable.
Any polynomial family $(f_t)_{t\in \dd C}$ as above 
containing $f$ has
all of its members $f_t$ non-linearizable except for 
an exceptional polar set of parameters $t$.

2) For an arbitrary holomorphic 
germ $\varphi(z)=\cl O (z^2)$ and for 
almost all values $t\in \dd C$ except a polar set $E$,
we have that 
$$
f_t (z)=e^{2\pi i \a } z +z^2 +t \varphi (z)
$$
is not linearizable.

3) Let 
$$
f(z)=e^{2\pi i \a } z +\sum_{n\geq 2} f_n z^n
$$
be an arbitrary entire function.
Keeping all coefficients fixed except $f_2$, there is 
a polar set $E$ such that if $f_2 \in \dd C-E$, then 
$f$ is not linearizable.}

Also, we have the same type of results for rational 
functions.

\eno {Corollary 4.}{Let 
$$
\cl R_{\lambda , d}=\{ f\in \dd C (z) ; f(0)=0; f'(0)=\lambda ; d^0 f=d \}
$$
When $\a \in \dd R-\dd Q$ is not a Bruno number, 
except for an exceptional set, all rational functions 
in $\cl R_{\lambda , d}$ are not linearizable.}

\medskip

The corollaries presented here are by no means restricted to 
dimension 1. Just one example of new result.

\eno {Corollary 5.}{ We consider the space 
$\cl P_{A, d}$ of polynomial germs of holomorphic diffeomorphisms
with non-resonant linear part $A$ of total degree $d$. The existence 
of one non-linearizable example forces all the others except an 
exceptional set of $\Gamma$-capacity $0$ to be non-linearizable.
This happens for instance when one eigenvalue of $A$ does not
satisfy Bruno's condition.}

We can prove also a version of theorem 1 for resonant linear
parts $A$ which has an independent interest (for example 
when applied to symplectic holomorphic mappings). When the 
linear part is resonant, the linearization is not uniquely 
determined. Nevertheless, given a polynomial family $(f_t)$ 
as in theorem 1 whose elements are all formally linearizable,
there always exist a canonical family of linearizations $(h_t)$ whose
coefficients depend polynomially on $t$ (see [PM3]).
The complete treatment of this situation requires some  
algebraic preliminaries. We do not  develop 
them in this article. We refer to [PM3] for a complete 
treatment. We content to prove here the following:

\eno {Theorem 2.}{We consider a family $(f_t)_{t\in \dd C^m}$
as in theorem 1 but we allow $A \in GL_n (\dd C )$ to be resonant.
We are also given a family of formal linearizations $(h_t)_{t\in \dd C^m}$
whose coefficients depend polynomially on $t$. We assume that the 
monomial of valuation $l$ has as coefficient a polynomial 
of degree bounded above by $C_0 +C_1 l$ for some $C_0 , C_1 >0$.

We have the following dichotomy:

1) The family $(f_t)_{t\in \dd C^m}$ is holomorphically linearizable
by the family $(h_t)_{t\in \dd C^m}$.

2) For all $t\in \dd C^m$ except for an exceptional
set $E$ of $\Gamma$-capacity $0$, $h_t$ is diverging.
}

One can also prove a statement similar to theorem 2 when 
$(f_t)$ is not formally linearizable but the family $(h_t)$ 
conjugates the family to a formal normal form ([PM3]). 
A particular relevant case is the one of a symplectic holomorphic
diffeomorphism with an elliptic fixed point. The formal 
conjugacy to Birkhoff's normal form is then in general diverging (see
[Si-Mo] section 30).
The formal normal form situation is also relevant 
when $A$ is not invertible.

\stit {Central manifolds.}

In situations where the dynamics is not linearizable, 
one can still have invariant manifolds through the 
fixed point (see for example [Pos], and [St] for a general
treatment in the case of holomorphic vector fields). 
Usually one has a formal equation whose coefficients
depend polynomially on the coefficients of $f_t$ thus 
on $t$. In these situations the following theorem 
applies.

\eno {Theorem 3.}{Under the same assumptions as
in theorem 1, we assume the existence of a formal 
invariant submanifold through $0$ with equation
$$
F_t (z)=0
$$
with $F_t : \dd C^n \to\dd C^p$ a formal mapping
whose coefficients depend polynomially on $t\in \dd C^m$.
More precisely, the coefficient of the monomial of valuation 
$l$ is a polynomial on $t$ of degree less than $C_0 +C_1 l$ where 
$C_0 , C_1 >0$ are constants.

We have the dichotomy:

(1) $F_t$ converges and defines an invariant submanifold
for all $t\in \dd C^m$.

(2) Except for an exceptional set of $\Gamma$-capacity $0$   
of parameter values $t\in \dd C^m$, $F_t$ diverges.
}

We have the same theorem for holomorphic vector fields.
To be more specific, consider the situation treated by L.Stolovitch
[Sto], for $1\leq j \leq n$,
$$
\dot z_j =\lambda_j z_j+\sum_{i=1}^d t^i f_{j,i} (z) 
$$
where $f_{j,i} =\cl O (2)$.
We assume that the linear part (which does not depend on $t$) is in 
the Siegel domain, that is $0$ belongs to the convex hull of 
$\{\lambda_1 ,\ldots , \lambda_n \}$. We assume that the linear 
part is resonant, and the resonances, $n_1 , \ldots , n_2 \geq 1$ and 
any $1\leq j \leq n$,
$$
\sum_{i=1}^n n_i \lambda_i -\lambda_j \not= 0 \ .
$$
are generated by a 
finite number of resonances, $1\leq j\leq l$, 
$r_j=(r_1, \ldots r_n)\not= 0$, $r_j\in \dd N^n$,
$$
(r_j ,\lambda) =0 \ .
$$
Then there exists a formal change of variables $w=h_t(z)$ with 
$h_t(0)=0$ and $D_0 h_t=I$ which transforms the system 
into 
$$
\dot w_i=\lambda_i w_i +g_{i,t} (w)
$$
with $g_{i,t} (w)=\sum_{j=1}^l g_{i,j,t} y^{r_j}$, and 
if $||r_j||=1$ then $g_{i,j,t} (0)=0$. As constructed in 
[Sto], the coefficients of the formal normalization do 
depend polynomially on $t$.

\eno {Theorem 4.}{ With the previous assumption, we have 
the following dichotomy,

\medskip

1) For all value of $t \in  \dd C^m$ the formal normalization 
$h_t$ converges, thus the sub-manifold $ \{ w^{r_1} =0 , \ldots 
w^{r_n}=0 \}$ is invariant.

\medskip

2) Except for an exceptional set of values of $t$ of $\Gamma$-capacity 
$0$, the normalization mappings $h_t$ diverge.
}

According to  [Sto], and assuming that the higher dimensional 
resonant Bruno condition on $(\lambda_1 , \ldots , \lambda_n )$
holds, we are always
in case (1).

\stit {Singularities of holomorphic vector fields.}

We consider a polynomial family of germs of holomorphic vector fields 
as before. But 
we assume here that the linear 
part is non-resonant, that is, for any $n_1 , \ldots , n_2 \geq 1$ and 
any $1\leq j \leq n$,
$$
\sum_{i=1}^n n_i \lambda_i -\lambda_j \not= 0 \ .
$$

\eno {Theorem 5.}{Under the above hypothesis, we have the 
dichotomy

\medskip

1) The family of holomorphic vector fields is linearizable for all $t$.

\medskip

2) Except for an exceptional set of values of $t$ of $\Gamma$-capacity 
$0$, the holomorphic vector fields are non-linearizable.
}

In the case $n=2$ one has a complete correspondence of the problem of 
linearization of holomorphic vector fields as above and the problem 
of linearization of germs of holomorphic diffeomorphisms of $(\dd C , 0)$
(see [PM-Yo] and the references there in).
Yoccoz and the author proved that Bruno condition is optimal for the 
problem of linearization.

\stit {Centralizers.}

We discuss here the situation of one complex variable. 
The analysis generalizes similarly to higher dimension.

The study of centralizers of holomorphic germs generalizes
the problem of linearization. We refer to [PM2] for proofs 
and references.
In the group of holomorphic
diffeomorphisms $G=({\hbox {\rm Diff}} (\dd C , 0), \circ )$,
composed by holomorphic germs $f$ with $f(0)=0$ and 
$f'(0)\not= 0$, we consider the centralizer of $f$,
$$
{\hbox {\rm Cent }} (f)=\{ g\in {\hbox {\rm Diff}} (\dd C , 0)
; g\circ f = f\circ g \}
$$
This group can be interpreted as the group of symmetries of 
$f$ (i.e. those changes of variables conjugating $f$ to 
itself). We have the following cases:

\medskip

\item {1)}For germs with attracting or repelling fixed point
at $0$, i.e.  $f'(0)=e^{2\pi i \a }$ with $\a \notin \dd R$,
the centralizer is a complex flow of dimension $1$.

\medskip

\item {2)}For germs with indifferent rational fixed point
at $0$, i.e.  $f'(0)=e^{2\pi i \a }$ with $\a \in \dd Q$,
the centralizer is generated by root (for composition) of 
the germ (then it is discrete), or it is a one dimensional 
complex flow.

\medskip

These cases are well understood.
We discuss the last case in what follows.

\medskip

\item {3)} For germs with an indifferent irrational fixed point 
at $0$, $f'(0)=e^{2\pi i \a }$ with $\a \in \dd R-\dd Q$, 
the centralizer can be a one-dimensional real flow (the linearizable 
case), discrete or uncountable. 
The occurrence of the last 
possibility was only proved recently in [PM2].

\medskip

In this case centralizer is abelian and isomorphic to 
a subgroup of the circle $\dd T =
\dd R /\dd Z$ by the rotation number morphism,
$$
\apl {\rho }{G} {\dd T} {f} {\log f'(0)}
$$
We denote 
$$
G(f) =\rho ({\hbox {\rm Cent}} (f) ) \ .
$$
Note that $\dd Z \a \subset G(f)$.
The holomorphic germ $f$ is linearizable if and only if 
the centralizer is full $G(f)=\dd T$, otherwise it is an
$F_{\sigma }$ and dense set of $\dd T$ with $0$ measure 
(and indeed $0$ capacity). Moreover, all elements 
$g\in {\hbox {\rm Cent }} (f)$ are non-linearizable.
 
Thus how small is $G(f)$ can be thought as a measure of 
how far is $f$ from being linearizable. Thus the study 
of centralizers (apart from the motivation coming 
from the theory of foliations, see [PM2]) is motivated
as a finer study of linearization.
The question of determining if $\beta \in G(f)$ 
is intimately connected with the common rational 
approximations of $\a$ and $\b$, as the following theorem 
of J. Moser shows ([Mo]). Let $f$ be non-linearizable.
If there exists $\g , \tau >0$ 
such that for any $p\geq 1$, $q \in \dd Z$,
$$
\min (|q\a -p_1| , |q\b -p_2| ) \geq {\gamma \over q^\tau }
$$
then $\beta \notin G (f)$.
The necessity of an arithmetic condition in Moser's theorem 
is proved in [PM2].

We have:

\eno {Theorem 6.}{Let $f_t$ be a family of holomorphic
germs as in theorem 1, with fixed linear part 
$f'(0)=e^{2\pi i \a }$, $\a \in \dd R-\dd Q$. 
For any $\b \in \dd T$, we have the following dichotomy:

1) For all $t\in \dd C$, $\b \in G(f_t)$.

2) Except for an exceptional polar set $E\subset  \dd C$,
$\b \notin G(f_t)$.
}

{\bf Further applications.}

A complete treatment for the problem of linearization of 
resonant holomorphic germs is given in [PM3]. These techniques
also apply to analytic K.A.M. of persistence of invariant tori.
In [PM4] we study the Lindsted series for the standard map. 
Behind the technique used here there is 
an abstract theorem on holomorphic extension of Rothstein type
for a certain type of power series. We discuss it in [PM4].




\null
\vfill
\eject


\stit {1) Proof of theorem 1.}

\stit {a) Nadirashvili and Bernstein lemmata.}

For the definition of Green function, polar sets and 
other notions in potential theory we refer the reader
to [Ra] for example (for a more encyclopedic treatment
see [Tsu]).

Y. Ilyashenko in his article [Il] makes use
of the following lemma attributed to N.S. Nadirashvili
([Na]). 

\eno {Lemma (Nadirashvili).}{Let $E\subset \dd C$ be a compact set 
with positive measure in the disk $\dd D_R$ of center $0$
and radius $R>0$. Let $P$ be a polynomial of degree $n$
such that for some $M>0$
$$
||P||_{C^0 (E)}\leq M^n \ .
$$
Then there exists a constant $C$ 
only depending on the measure of $E$ and $R>0$ such that 
$$
||P||_{\dd D_R } \leq C^n M^n \ .
$$
}

We improve on [Il]  
observing that 
the measure of $E$ is not the relevant quantity.
Nadirashvili's lemma is a direct corollary of the classical Bernstein
lemma in approximation theory and classical potential theory
(see [Ra] p.156) and   
the fact that a set of positive measure
is non-polar :

\eno {Lemma (Bernstein).}{Let $E\subset \dd C$ be a non-polar
compact set (i.e. ${\hbox {\rm cap}} (E) >0$). Let $\O$ be the connected 
component of ${\overline {\dd C}}-E$ containing $\infty$.
Then for any polynomial $P$ of degree $n$, we have for $t\in \dd C$,
$$
|P(t)|\leq e^{n g_{\O} (t,\infty)} \ ||P||_{C^0(K)} 
$$
where $g_{\O}$ denotes the Green function of $\O$.
}

The proof is quite simple.

\dem {Proof.}{We can assume the polynomial monic. Then
$$
u(t)=\log P(t)-\log ||P||_{C^0(K)} -g_{\O } (t, \infty)
$$
is sub-harmonic, is negative near $\infty$ (because 
$g_{\O }(t, \infty )=\log |t| +{\hbox {\rm cap}} (E) +o(1)$),
and $\limsup u(t) \leq 0$ when $t\to K$. The maximum principle
concludes the proof.}

For future reference we recall here that a countable 
union of polar sets is polar.

\stit {b) $\Gamma$-capacity.}

We recall the definition of $\Gamma$-capacity and 
we refer to [Ro] for more properties.
Let $E \subset \dd C^m$. The $\Gamma$-projection of $E$ on 
$\dd C^{m-1}$ is the set $\Gamma_m^{m-1} (E)$ 
of $z=(z_1, \ldots , z_{m-1})\in \dd C^{m-1}$
such that 
$$
E\cap \{ (z, w)\in \dd C^m \}
$$
has positive capacity in the complex plane 
$ \dd C_z=\{ (z, w)\in \dd C^m \}$.
We define 
$$
\Gamma_m^1 (E)=\Gamma_2^1 \circ \Gamma_3^2\circ \ldots \Gamma_m^{m-1} (E) \ .
$$
Finally, the $\Gamma$-capacity is defined as 
$$
\Gamma {\hbox {\rm -Cap}} (E)=\sup_{A\in U(m,\dd C)} {\hbox {\rm Cap}}
\ \Gamma_m^1 (A(E)) \ .
$$
where $A$ runs over all unitary transformations of $\dd C^m$.

Using the definition of $\Gamma$-capacity it is easy to see 
that we are reduced to prove the theorems for $m=1$

\stit {c) Proof of theorem 1.}

We have the following elementary lemma.

\eno {Lemma 1.1.}{The coefficient vectors $h_i(t)$ of the
formal linearization 
$$
h_t(z)=z+\sum_{i=(i_1, \ldots i_n)\atop 
i_1+\ldots +i_n \geq 2} h_i(t) z^i
$$
have coordinates that are polynomials in the parameter 
$t=(t_1, \ldots t_m)$ of degree less than $d (i_1+\ldots i_n)$.}

\dem {Proof.}{We can assume that $A$ is in 
upper triangular Jordan normal form.
We solve the functional equation
$$
A\circ h_t=h_t\circ f_t 
$$
identifying coordinates and developing 
in homogeneous vector monomials. By induction on 
$|i|=i_1+\ldots +i_n$ we do determine successively 
the vectors $h_i(t)$ that depend on coefficients of $f_t$
and lower order $h_j(t)$'s, $|j|<|i|$. By induction,
the linear equations
determining $h_i(t)$ do have the form
$$
(A-M_i) h_i(t)= \sum_{|j|<|i|} c_j(t) h_j(t)
$$
where the matrix $M_i$ is upper triangular, 
only depends on $A$ and $i$ (but 
not on the parameter $t$), has diagonal coefficients
products of eigenvalues of $A$ (thus $A-M_i$ is 
invertible) and in the left hand side the 
coefficients $c_j(t)$  are polynomials in $t$ of total 
degree at most $d$. By induction the result follows.}

\dem {Proof of theorem 1.}{According to previous section 
we are reduced to prove the theorem to the case $m=1$.
Let 
$$
E=\{ t\in \dd C ; f_t \ {\hbox {\rm    is linearizable  }} \} \ .
$$
We want to show that $E$ is polar or the whole complex plane.
We have 
$$
E=\bigcup_{j\geq 1} E_j
$$
where $E_j$ the set of parameters $t$ such that $h_t$ has radius
of convergence larger or equal to $1/j$. Thus if $E$ is non-polar, 
we have that for some 
$j\geq 1$, $E_j$ is not polar. 
Thus there exists $\rho_0 >0$
such that for all $t\in E_j$, 
$$
\varphi (t)=\limsup_{|i|\to +\infty} ||h_i(t)||\rho_0^{-|i|} <+\infty \ .
$$
The function $\varphi$ is lower semicontinuous, and
$$
E_j=\bigcup_p L_p
$$
where $L_p=\{ z\in E_j ; \varphi (t) \leq p$ is closed.  
By Baire theorem for some $p$, $L_p$ has non-empty interior (with 
respect to $E_j$), thus this $L_p$ has positive capacity. Finally
we found a compact set $C=L_p$ of positive capacity such that 
there exists $\rho_1 >0$ such that for any $t\in C$ and 
and all $i\in \dd N^n$,
$$
||h_i (t)||_{C^0 (C)} \leq \rho_1^{|i|} \ .
$$
Using Bernstein's lemma and lemma 1.1 we get that for any compact 
set $K\subset \dd C$ we have
$$
||h_i(t)||_{C^0 (K)} \leq C(K)^{d |i|} \rho_1^{|i|} \ ,
$$
for some constant $C(K)$ depending only on $K$.
Thus $f_t$ is linearizable for any $t\in \dd C$.
The constant $C(K)$ can be estimated by the precise form 
of Bernstein lemma as
$$
C(K) = \exp (\sup_K g_{\O} (t , \infty )
$$
where $\O$ is the connected component containing $\infty$ of
the complement of $C$. The asymptotic
$$
g_\O (t, \infty) \approx \log |t|
$$
for $t\to \infty$ gives the lower estimate on the 
radius of convergence.
}

\stit {c) Proof of the corollaries.}

Corollary 1 is just a particular case of theorem 1. Corollary 3
part (1) also. Now using corollary 3 part (1) we prove

\eno {Corollary (Yoccoz, [Yo]).}{
The quadratic polynomial $P_{\a} (z)=e^{2\pi i \a } z +z^2$
is non-linearizable when $\a$ is not a Bruno
number.}

\dem {Proof.} {For this, pick $f$ non-linearizable (that exists from [Yo])
and consider $f_t(z)=t P_{\a} (z) +(1-t)f(z)$. Since $f_1$ is not 
linearizable, all $f_t$ except for a polar set $E$ of values of $t$
are not linearizable. By Douady-Hubbard straitening theorem
$\dd C-E$ is a neighborhood of $0$ and $0$ is not linearizable.}

Now by the same argument, part (2) and (3) of corollary 3 follow.

Now we prove:

\eno {Corollary (P\'erez-Marco, [PM1]).}{ If $P$ is a 
structurally stable polynomial in the space
$$
\cl P_{\lambda , d}=\{ P(z)=\lambda z+a_2 z^2 +\ldots +a_d z^d ; 
(a_2, \ldots , a_d ) \in \dd C^{d-1} \} 
$$
then $P$ is not linearizable.}

\dem {Proof.}{Just consider 
$$
f_t (z)=t P(z)+(1-t) P_{\a} (z)
$$
and do the same proof.}

Now corollary 2 follows from part (1) of corollary 3. 

\medskip

Corollary 4 is not a strict corollary of theorem 1 but 
exactly the same proof applies, observing that the 
coefficients of the linearization are polynomial functions
of the coefficients of the rational function with appropriate
degree.

\medskip

For corollary 5 only the last assertion is not immediate.
If one of the eigenvalues of $A$ violates Bruno condition, 
$\lambda_1$ for example, then 
$$
(z_1 , \ldots , z_n)\mapsto (\lambda_1 z_1+z_1^2 , \lambda_2 z_2 , 
\ldots , \lambda_n z_n )
$$
is not linearizable, thus the first part applies giving a rich
family of polynomial non-linearizable examples.

\stit {2) Proof of the other theorems.}

\stit {a) Formal linearizations and theorem 2.}

The proof of theorem 2 is similar to the proof given in 
the previous section of theorem 1. For a complete study 
of the resonant case we refer to [PM3]. We just mention 
here the new difficulties that appear.

Assume that we have a germ of holomorphic diffeomorphism
$$
f: (\dd C^n , 0) \to (\dd C^n,0) 
$$
with resonant linear part $A=D_0f \in GL_n(\dd C )$.

The formal linearization $h$ is not always unique 
when the linear part $A$ is resonant or not invertible.
For example, for $n=2$ and 
$$
A=\pmatrix { 1 &1 \cr 0 & 1\cr}
$$
then if $h$ is a formal linearization then $l\circ h$ is 
also one, where
$$
l(z_1, z_2)= (z_1+k(z_2) , z_2)\ .
$$
One can show that all linearizations can be obtained 
that way ([PM3]).
If we consider a polynomial family $(f_t)_{t\in \dd C^m}$, 
to request that $h_t$ has coefficients depending polynomially 
on $t$ does not improve things. One then can take various
$k_t$ depending polynomially on $t$. Thus the family $(h_t)$
with this further restriction is not unique. 

This presents a subtle problem in order to prove the non-linearizability.
Considering a polynomial parameter family of formal linearizations of 
a fixed map $f$, $(h_t)$, we may be in the second case, but this does
not mean that $f$ is not linearizable. For instance, if the 
exceptional set $E$ is not empty, then $f$ will be linearizable !
The question of non-linearizability is harder 
to answer. In [PM3] we show that if the polynomial family of linearizations
is chosen in a natural way, this difficulty does not arise.

\stit {b) Other theorems.}

The proofs are similar than in section 1. We just comment on the 
particularities of each problem.

For an explicit example where theorem 3 applies one
can workout the example of J. Poschel [Pos]. The polynomial 
dependence with the appropriate bound on the 
degrees follows from the formal computation of the 
formal equation of the invariant manifold. 

\medskip

Theorem 4 is proved in a similar way. We refer to [Sto] for 
the formal computation of a normalizing map with polynomial 
dependence on the parameter $t$ with the appropriate degrees. 
One can workout in this situation 
similar results than in [PM3].
 
\medskip

The linearization in theorem 5 is unique and it is well known ([Ar]) that 
it depends polynomially on $t$ with the appropriate degrees. Thus the same
proof applies. Note that in $\dd C^2$, by [PM-Y] one can realize any germ of holomorphic
diffeomorphism in $(\dd C , 0)$ as holomony of a singularity of 
holomorphic vector field of the type considered. The realization of 
the quadratic map is also structurally stable in the following sense: 
Any nearby holomorphic vector field has a holonomy that is quasi-conformally
conjugated to the quadratic polynomial. Thus for any one parameter polynomial 
family containing this vector field, if $\a $ is not a Bruno number, we 
will be in case 2. 

\medskip

For the proof of theorem 6, we give the induction formulas for 
the coefficients of $g_\b $. Let $\mu =e^{2\pi i \a }=f_1$ and
$\lambda=e^{2\pi i \a }=g_1$, and
$$
\eqalign {f(z) &= \sum_{n=1}^{+\infty } f_n z^n \cr
g(z) &= \sum_{n=1}^{+\infty } g_n z^n \cr
}
$$
Identifying terms of degree $n\geq 2$
in the equation $g\circ f =f\circ g$, we get for $n\geq 2$,
$$
g_n ={\mu^n -\mu \over \lambda^n -\lambda} f_n +
\sum_{p=1}^{+\infty } f_p \sum_{i_1+
\ldots i_p =n \atop i_j \geq 1} g_{i_1}
\ldots g_{i_p}  +
\sum_{p=1}^{+\infty } g_p \sum_{i_1+\ldots i_p =n
\atop i_j \geq 1} f_{i_1}
\ldots f_{i_p}  
$$
And by induction the coefficients of $g_\b$ depend 
polynomially on $t$ and have the appropriate degrees.


\bigskip

{\centerline {\bf BIBLIOGRAPHY}}

\bigskip
\bigskip

\ref{[Ar]}{V.I. ARNOLD}{Chapitres supplementaires de la th\'eorie
des \'equations diff\'erentielles ordinaires}{Editions MIR, 
Moscou, 1980.}

\ref{[Cr]}{H. CREMER}{\"Uber die H\"aufigkeit der 
Nichtzentren}{Math. Ann. {\bf 115} (1938),
p. 573-580.}

\ref{[He1]}{M. R. HERMAN}{Une m\'ethode pour minorer les
exposants de Lyapounov et quelques exemples montrant
le caract\`ere local d'un th\'eor\`eme d'Arnol'd et Moser
sur le tore en dimension 2}{Comment. Math. Helv. {\bf 58}, 
1983, p.453-502.}

\ref{[He2]}{M. R. HERMAN}{Recent results and some 
open questions on Siegel's lin\'earisation theorem 
of germs of complex analytic diffeomorphisms 
of $\dd C^n$ near a fixed point}{Proceedings 
$VIII^{th}$ Int. Conf. Math. Phys., World Scientific 
Publishers, Singapour, 1987.}

\ref{[Il]}{Y. ILYASHENKO}{Divergence of series that reduce an analytic 
differential equation to linear normal form at a singular point}{Functional 
Analysis and Applications, {\bf 13}, 1979, p.87-88.}

\ref{[Na]}{N.S. NADIRASHVILI}{}{Vest. Mosk. Gos. Univ., 
Ser. Mat., {\bf 3}, 1976, p.39-42. From [Il]).}

\ref{[Poi]}{H. POINCAR\'E}{Les m\'ethodes nouvelles de
la m\'ecanique c\'eleste}{Volume 1, Paris, 1893.}

\ref{[Pos]}{J. POSCHEL}{On invariant manifolds of complex analytic 
mappings near fixed points}{Exposition. Math. {\bf 4}, 2, 1986, p.97-109.}

\ref{[PM1]}{R. P\'EREZ MARCO}{Sur les dynamiques holomorphes 
non lin\'earisable et une conjecture de V. I. Arnold} 
{Ann. Scient. Ec. Norm. Sup. 4 serie, t.26, 1993, p.565-644.}

\ref{[PM2]}{R. P\'EREZ MARCO}{Non linearizable holomorphic dynamics
having an uncountable number of symmetries}{Invent. Math., {\bf 119},
 1995, p. 67-127.}

\ref{[PM3]}{R. P\'EREZ MARCO}{Linearization of holomorphic germs with 
resonant linear part} 
{preprint, 2000.}

\ref{[PM4]}{R. P\'EREZ MARCO}{Total convergence or
general divergence of the Lindsted series of the standard map} 
{preprint, 2000.}

\ref{[PM5]}{R. P\'EREZ MARCO}{A note on holomorphic extensions} 
{preprint, 2000.}

\ref{[PM-Yo]}{J.-C. YOCCOZ, R. P\'EREZ MARCO}{Germes de feuilletages 
holomorphes \`a holonomie prescrite}{"Complex methods in dynamical systems" 
Ast\'erisque {\bf 222}, 1994, p.345-371.}

\ref{[Ra]}{T. RANSFORD}{Potential theory in the complex plane}{London
Mathematical Society, Student Texts {\bf 28}, Cambridge University 
Press, 1995.}

\ref{[Si-Mo]}{K.L. SIEGEL, J. MOSER}{Lectures on celestial mechanics}
{Springer-Verlag, {\bf 187}, 1971.}

\ref{[Sto]}{L. STOLOVITCH}{Sur un th\'eor\`eme de Dulac}
{Ann. Inst. Fourier, {\bf 44}, 5, 1994, p.1397-1433.}

\ref{[Tsu]} {M. TSUJI}{Potential Theory in Modern Function 
Theory}{Maruzen, Tokyo, 1959.}

\ref{[Yo]}{J.-C. YOCCOZ}{Th\'eor\`eme de Siegel, 
polyn\^omes quadratiques et nombres de Brjuno}
{Ast\'erisque {\bf 231}, 1995.}

\end